\documentclass[10pt]{article}
\usepackage{amssymb}
\usepackage{amsmath,amssymb}\usepackage{cite}\usepackage{mathrsfs}\usepackage[amsmath,thmmarks]{ntheorem}
\usepackage{listings}
\usepackage[titletoc]{appendix}\allowdisplaybreaks

\makeatletter
\renewcommand\theequation{\thesection.\@arabic\c@equation}
\makeatother
\newtheorem{thm}{ Theorem}[section]%
\newtheorem{lem}[thm]{ Lemma}%
\newtheorem{Pro}[thm]{ Proposition}%
\newtheorem{Que}[thm]{Question}%
\newtheorem{Con}[thm]{Conjecture}%
\newtheorem{Fac}[thm]{Fact}%
\input cyracc.def

\def\f{\noindent}

\def\demo{\f{\bf Proof}\hskip10pt}

\def\qed{\hfill $\Box$}

\begin{document}
\title{\textbf{Automorphism group of the graph $A(n,k,r)$\\{{\small Dedicated to my PhD supervisor Xiuyun Guo's 70th birthday.}}}}

\footnotetext{
The work was supported by Hainan Provincial Natural Science Foundation of China (No. 122RC652).
\\ \ \ \ \ E-mail address: Junyao$_{-}$Pan@126.com.}

\author{Junyao Pan\,
 \\\\
School of Sciences, Wuxi University, Wuxi, Jiangsu,\\ 214105
 People's Republic of China \\}
\date {} \maketitle

\baselineskip=16pt

\vskip0.5cm

{\bf Abstract:} Let $[n]^{(k)}$ be the set of all ordered $k$-tuples of distinct elements in $[n]=\{1,2,...,n\}$. The $(n,k,r)$-arrangement graph $A(n,k,r)$ with $1\leq r\leq k\leq n$, is the graph with vertex set $[n]^{(k)}$ and with two $k$-tuples are adjacent if they differ in exactly $r$ coordinates. In this manuscript, we characterize the full automorphism groups of $A(n,k,r)$ in the cases that $1\leq r=k\leq n$ and $r=2<k=n$. Thus, we resolve two special cases of an open problem proposed by Fu-Gang Yin, Yan-Quan Feng, Jin-Xin Zhou and Yu-Hong Guo. In addition, we conclude with a bold conjecture.



\section {Introduction}
Throughout this paper, $[n]$ denotes the standard $n$-element set $\{1,2,...,n\}$; $[n]^{(k)}$ stands for the set of all ordered $k$-tuples of distinct elements in $[n]$; $S_n$ expresses the symmetric group on $[n]$. In addition, some notions and notations about permutation groups can be seen in \cite{C,D}.

The \emph{$(n,k)$-arrangement graph} $A(n,k)$ with $1\leq k\leq n-1$, is defined to have vertex set $[n]^{(k)}$ with two vertices $[s_1,s_2,...,s_k]$ and $[t_1,t_2,...,t_k]$ are adjacent if and only if $|\{i|s_i\neq t_i,1\leq i\leq k\}|=1$. Actually, the arrangement graph was firstly introduced and studied by Day and Tripathi in \cite{DT} as interconnection network topology, which generalizes the star graph topology in \cite{AK}. Since then a fair amount of work has been done on the arrangement graphs in the literature, for example \cite{CLY,CGW,HLTH,YFZG}. Chen et al \cite{CGW} generalized the $(n,k)$-arrangement graph $A(n,k)$ to so called \emph{$(n,k,r)$-arrangement graph} $A(n,k,r)$. The $(n,k,r)$-arrangement graph $A(n,k,r)$ with $1\leq r\leq k\leq n$ has vertex set $[n]^{(k)}$ with two vertices $[s_1,s_2,...,s_k]$ and $[t_1,t_2,...,t_k]$ are adjacent if and only if $|\{i|s_i\neq t_i,1\leq i\leq k\}|=r$. Yin et al \cite{YFZG} not only investigated the symmetry property of $(n,k)$-arrangement graph but also proposed the following open problem.

\begin{Que}\label{pan1-1}\normalfont(\cite[Question~5.1]{YFZG}~)
Determine the automorphism group of the $(n,k,r)$-arrangement graph $A(n,k,r)$ with $2\leq r\leq k\leq n$.
\end{Que}

Given an undirected graph $\Gamma = (V, E)$ with vertex set $V$ and edge set $E$. Let $f$ be a bijection from $V$ to $V$ such that $(\alpha^f,\beta^f)\in E$ if and only if $(\alpha,\beta)\in E$ for all $\alpha,\beta\in V$. Then $f$ is called an \emph{automorphism} of $\Gamma$, and the \emph{full automorphism group} of $\Gamma$ is denote by $Aut(\Gamma)$. Actually, the research on the automorphism groups of some graphs has always been an interesting topic, for examples \cite{F,FMW,G,GR,KNZ,P,PG,T}. Thereby, we are interested in Question\ \ref{pan1-1}. Consider $A(n,k,r)$ with $2\leq r\leq k\leq n$. Motivated by \cite{DT,YFZG}, for each $g\in S_n$, define the permutation $P(g)$ on $[n]^{(k)}$ by $$[s_1,s_2,...,s_k]^{P(g)}=[s^g_1,s^g_2,...,s^g_k],[s_1,s_2,...,s_k]\in[n]^{(k)};$$ for each $h\in S_k$, define the permutation $Q(h)$ on $[n]^{(k)}$ by $$[s_1,s_2,...,s_k]^{Q(h)}=[s_{1^{h^{-1}}},s_{2^{h^{-1}}},...,s_{k^{h^{-1}}}],[s_1,s_2,...,s_k]\in[n]^{(k)}.$$  Clearly, $P(S_n)$ and $Q(S_k)$ are also two subgroups of $Aut(A(n,k,r))$. Furthermore, one easily checks that $P(S_n)Q(S_k)=P(S_n)\times Q(S_k)$ and $B_{[1,2,...,k]}\cong S_k\times S_{n-k}$, where $B_{[1,2,...,k]}$ is the stabilizer of $P(S_n)Q(S_k)$ at $[1,2,...,k]$. In other words, the \cite[Lemma~3.1,~3.2]{YFZG} still hold when $1<r\leq k\leq n$. However, it seems to be rather difficult to obtain some results as the \cite[Lemma~3.4]{YFZG} which plays a crucial role in proving \cite[Theorem~1.1]{YFZG}. Perhaps for this reason, there has been no progress on this open problem so far.

In this manuscript, we investigate the Question\ \ref{pan1-1} from another perspective. Let $\Gamma = (V, E)$ be an undirected graph with vertex set $V$ and edge set $E$. If $S\subseteq V$ such that there does not exist edge between $\alpha$ and $\beta$ for any two vertices $\alpha,\beta\in S$, then $S$ is called an \emph{independent set} of $\Gamma$. If $S$ is an independent set while $S\cup\{\gamma\}$ is not an independent set for each $\gamma\in V\setminus S$, then $S$ is called a \emph{maximal independent set}. If $S$ is a maximal independent set such that $|S|\geq|T|$ for any other maximal independent set $T$, then $S$ is called a \emph{maximum independent set}. Let $\Omega$ denote the set of all maximum independent sets of $\Gamma$. Obviously, $Aut(\Gamma)$ induces a permutation group on $\Omega$. This is precisely our breakthrough point, which we resolve Question\ \ref{pan1-1} in the cases that $1\leq r=k\leq n$ and $r=2<k=n$ by constructing maximum independent sets and isomorphic graphs.

\begin{thm}\label{pan1-2}\normalfont
If $n>2$, then
\begin{align*}
Aut(A(n,k,r))&\cong
  \begin{cases}
    P(S_n)\times Q(S_k) & \text{if } r=k<n, \\
    [P(S_n)\times Q(S_n)]\rtimes Z_2 & \text{if } r=k=n, \\
    [P(S_n)\times Q(S_n)]\rtimes Z_2 & \text{if } r=2,k=n.
  \end{cases}
\end{align*}
Here, $Z_2=\langle h\rangle$ and $h$ is the map that $[s_1,s_2,...,s_n]^h=[x_1,x_2,...,x_n]$ where $x_j=i$ if $s_i=j$ for $j=1,2,...,n$.
\end{thm}

However, the maximum independent set of $A(n,k,r)$ becomes complicated when $r< k$. So we cannot determine the automorphism group of $A(n,k,r)$ in this case except $r=2,k=n$. But, in the light of some special phenomena, we conclude with some remarks and a bold conjecture.

\section {Proof of Theorem\ \ref{pan1-2}}

 Now we start from constructing maximum independent sets of $A(n,k,k)$. For a pair of number $(i,j)$ with $i\in[n]$ and $j\in[k]$, we define a subset of $[n]^{(k)}$ as follows:$$\Delta_{ij}=\{[s_1,s_2,...,s_j,...,s_k]|s_j=i,s_t\in[n]\setminus\{i\}~{\rm{for~all}}~t\in[k]\setminus\{j\}\}.$$
For convenience, we set $$\Omega=\{\Delta_{ij}|i\in[n],j\in[k]\}.$$
Clearly, $\Delta_{ij}$ is a maximal independent set of $A(n,k,k)$. Indeed, we observe that $\Delta_{ij}$ is a maximum independent set of $A(n,k,k)$, and there is the following result.

\begin{Pro}\label{pan2-1}\normalfont
Let $S\subseteq[n]^{(k)}$ with $n>2$. Then $S$ is a maximum independent set of $A(n,k,k)$ if and only if $S\in\Omega$.
\end{Pro}
\demo Let $[s_1,s_2,...,s_k]$ and $[t_1,t_2,...,t_k]$ be two distinct $k$-tuples in $[n]^{(k)}$. If there exists an $i\in[k]$ such that $s_i=t_i$, then we say that $[s_1,s_2,...,s_k]$ and $[t_1,t_2,...,t_k]$ are intersecting. Thus, $[s_1,s_2,...,s_k]$ and $[t_1,t_2,...,t_k]$ are not adjacent in $A(n,k,k)$ if and only if $[s_1,s_2,...,s_k]$ and $[t_1,t_2,...,t_k]$ are intersecting. Therefore, $S$ is an independent set of $A(n,k,k)$ if and only if $S$ is an intersecting family. On the other hand, $[s_1,s_2,...,s_k]$ can be viewed as an injection from $[k]$ to $[n]$ under the rule $i\mapsto s_i$ for $i=1,2,...,k$. By \cite{FS}, $[s_1,s_2,...,s_k]$ and $[t_1,t_2,...,t_k]$ are intersecting if and only if they are intersecting as injections.
Apply \cite[Theorem~2.1]{FS}, we deduce this proposition immediately.    \qed

It follows from Proposition\ \ref{pan2-1} that $Aut(A(n,k,k))$ induces a permutation group acting on $\Omega$, denoted by $Aut(A(n,k,k))|_\Omega$. It is well-known that there exists a natural surjective homomorphism $$\varphi:~Aut(A(n,k,k))\rightarrow Aut(A(n,k,k))|_\Omega.$$
Apply the fundamental homomorphism theorem, we deduce that $$Aut(A(n,k,k))/ker\varphi\cong Aut(A(n,k,k))|_\Omega,$$
where the kernel $Ker\varphi$ consists of all automorphisms in $Aut(A(n,k,k))$ which fix every $\Delta_{ij}$ in $\Omega$. In particular, if $ker\varphi=1$ then $Aut(A(n,k,k))\cong Aut(A(n,k,k))|_\Omega$, where $ker\varphi=1$ means that the kernel of $\varphi$ is trivial. In fact, we discover that $ker\varphi=1$ for all $n>2$. Next we recall some notions and notations which are used in proving this discovery.

Let $\alpha$ and $\beta$ be two vertices of an undirected graph $\Gamma = (V, E)$ with vertex set $V$ and edge set $E$. If there exists an edge between $\alpha$ and $\beta$ in $\Gamma $, then $\alpha$ and $\beta$ are called \emph{neighbor} (see \cite{R}). We use $N(\alpha)$ to stand for the collection of all neighbours of vertex $\alpha$. In addition, for $S\subseteq V$, we define $N(S)=\bigcap\limits_{\alpha\in S}N(\alpha)$. Next we start to prove $ker\varphi=1$ from a well-known fact.

\begin{Fac}\label{pan2-0}\normalfont
Let $\Gamma = (V, E)$ be an undirected graph. Then for any $S\subseteq V$ and $g\in Aut(\Gamma)$,
$$N(S)^g=N(S^g).$$
\end{Fac}

\begin{Pro}\label{pan2-2}\normalfont
Let $\varphi$ be the natural surjective homomorphism $$\varphi:~Aut(A(n,k,k))\rightarrow Aut(A(n,k,k))|_\Omega$$ with $n>2$. Then $ker\varphi=1$.
\end{Pro}
\demo Proof by contradiction. Suppose that there exists an automorphism $f\in Aut(A(n,k,k))$ such that $f$ is a non-identity while fixes every $\Delta_{ij}$ in $\Omega$. Assume that $[s_1,s_2,...,s_k]^f=[t_1,t_2,...,t_k]$ with $[s_1,s_2,...,s_k]\neq[t_1,t_2,...,t_k]$. Thereby, $s_i\neq t_i$ for some $i$. Note that $$N([s_1,s_2,...,s_k])\cap\Delta_{t_ii}\neq\emptyset~{\rm{and}}~N([t_1,t_2,...,t_k])\cap\Delta_{t_ii}=\emptyset.$$
On the other hand, by Fact\ \ref{pan2-0} we see that $N([s_1,s_2,...,s_k])^f=N([t_1,t_2,...,t_k])$. Therefore,
$$(N([s_1,s_2,...,s_k])\cap\Delta_{t_ii})^f=N([s_1,s_2,...,s_k]^f)\cap(\Delta_{t_ii})^f=N([t_1,t_2,...,t_k])\cap\Delta_{t_ii}.$$
This is a contradiction, as desired.    \qed

Let $G$ be a transitive permutation group acting on $[n]$. A non-empty subset $\Theta$ of $[n]$ is called a \emph{block} for $G$ if for each $g\in G$ either $\Theta^g=\Delta$ or $\Theta^g\cap\Theta=\emptyset$. Clearly, the singletons $\{i\}$ $(i\in[n])$ and $[n]$ are blocks, and these blocks are called \emph{trivial} blocks. Any others are called \emph{nontrivial} blocks. Put $\Sigma=\{\Theta^g|g\in G\}$, where $\Theta$ is a block of $G$. We call $\Sigma$ the \emph{system of blocks} containing $\Theta$.

\begin{Pro}\label{pan2-3}\normalfont
Let $\Sigma=\{\Omega_1,\Omega_2,...,\Omega_n\}$, where $\Omega_i=\{\Delta_{i1},\Delta_{i2},...,\Delta_{ik}\}$ for $i\in[n]$. If $k<n$, then $\Sigma$ is a system of blocks of $Aut(A(n,k,k))|_\Omega$ acting on $\Omega$.
\end{Pro}
\demo It suffices to prove that either $(\Omega_i)^g=\Omega_i$ or $(\Omega_i)^g\cap\Omega_i=\emptyset$ for every $g\in Aut(A(n,k,k))$ and $i\in[n]$. Proof by contradiction. Suppose that there exists an automorphism $g\in Aut(A(n,k,k))$ such that $(\Omega_i)^g\neq\Omega_i$ and $(\Omega_i)^g\cap\Omega_i\neq\emptyset$ for some $i$. In other words, there exist $a,b,c,d,j$ such that $(\Delta_{ia})^g=\Delta_{ib}$ and $(\Delta_{ic})^g=\Delta_{jd}$ and $i\neq j$ and $a\neq c$.

Since $a\neq c$, we deduce that $\Delta_{ia}\cap\Delta_{ic}=\emptyset$. Thereby, $\Delta_{ib}\cap\Delta_{jd}=\emptyset$, and so $b=d$. In addition, $k<n$ implies that we can take $[s_1,...,s_{b-1},i,s_{b+1},...,s_k]$ and $[s_1,...,s_{b-1},j,s_{b+1},...,s_k]$ in $\Delta_{ib}$ and $\Delta_{jd}$ respectively. Suppose $\alpha^g=[s_1,...,s_{b-1},i,s_{b+1},...,s_k]$ and $\beta^g=[s_1,...,s_{b-1},j,s_{b+1},...,s_k]$. Clearly, $\alpha\in\Delta_{ia}$ and $\beta\in\Delta_{ic}$. Thus, we can set $$\alpha=[t_1,...,t_{a-1},i,t_{a+1},...,t_k]~{\rm{ and}}~\beta=[r_1,...,r_{c-1},i,r_{c+1},...,r_k].$$ Note that $$[x_1,x_2,...,x_k]\in N([s_1,...s_{b-1},i,s_{b+1},...,s_k])\cap N([s_1,...s_{b-1},j,s_{b+1},...,s_k])$$ if and only if $x_b\not\in\{ i, j\}$ and $x_m\not\in\{s_m\}$ for each $m\in[k]\setminus\{b\}$. By the same token, $$[y_1,y_2,...,y_k]\in N([t_1,...,t_{a-1},i,t_{a+1},...,t_k])\cap N([r_1,...,r_{c-1},i,r_{c+1},...,r_k])$$ if and only if $y_a\not\in\{ i, r_a\}$ and $y_c\not\in\{ i, t_c\}$ and $y_m\not\in\{r_m,t_m\}$ for each $m\in[k]\setminus\{a,c\}$. Since $a\neq c$, we have $r_a\neq i$ and $s_c\neq i$. Thereby, we see that $$|[N(\alpha)\cap N(\beta)]^g|<| N([s_1,...s_{b-1},i,s_{b+1},...,s_k])\cap N([s_1,...s_{b-1},j,s_{b+1},...,s_k])|.$$
However, this is contradict to the Fact\ \ref{pan2-0}. Thus, $\Sigma$ is a system of blocks of $Aut(A(n,k,k))|_\Omega$ acting on $\Omega$.  \qed

\begin{Pro}\label{pan2-4}\normalfont
Let $\Sigma'=\{\Omega^1,\Omega^2,...,\Omega^k\}$, where $\Omega^j=\{\Delta_{1j},\Delta_{2j},...,\Delta_{nj}\}$ for $j\in[k]$. If $k<n$, then $\Sigma'$ is a system of blocks of $Aut(A(n,k,k))|_\Omega$ acting on $\Omega$.
\end{Pro}
\demo It suffices to prove that either $(\Omega^j)^g=\Omega^j$ or $(\Omega^j)^g\cap\Omega^j=\emptyset$ for every $g\in Aut(A(n,k,k))$ and $j\in[k]$. Proof by contradiction. Suppose that there exists an automorphism $g\in Aut(A(n,k,k))$ such that $(\Omega^j)^g\neq\Omega^j$ and $(\Omega^j)^g\cap\Omega^j\neq\emptyset$ for some $j$. Therefore, there exist $a,b,c,d,i$ such that $(\Delta_{aj})^g=\Delta_{bj}$ and $(\Delta_{cj})^g=\Delta_{di}$ and $i\neq j$ and $a\neq c$.

Since $a\neq c$, we deduce that $\Delta_{aj}\cap\Delta_{cj}=\emptyset$. Therefore, $\Delta_{bj}\cap\Delta_{di}=\emptyset$, and thus $b=d$. Moreover, $k<n$ indicates that we can pick $[s_1,...,s_{j-1},a,s_{j+1},...,s_k]$ and $[s_1,...,s_{j-1},c,s_{j+1},...,s_k]$ in $\Delta_{aj}$ and $\Delta_{cj}$ respectively. Suppose $[s_1,...,s_{j-1},a,s_{j+1},...,s_k]^g=[t_1,...,t_{j-1},b,t_{j+1},...,t_k]$ and $[s_1,...,s_{j-1},c,s_{j+1},...,s_k]^g=[r_1,...,r_{i-1},d,r_{i+1},...,r_k]$. Proceeding as in the proof of Proposition\ \ref{pan2-3}, we deduce that $\Sigma'$ is also a system of blocks of $Aut(A(n,k,k))|_\Omega$ acting on $\Omega$. \qed

\begin{lem}\label{pan2-5}\normalfont
If $n>k\geq1$, then $Aut(A(n,k,k))\cong P(S_n)\times Q(S_k)$.
\end{lem}
\demo Clearly, this lemma holds for $k=1$. Suppose that $n>k>1$. It follows from Proposition\ \ref{pan2-3} that $\Sigma$ is a system of blocks of $Aut(A(n,k,k))|_\Omega$ acting on $\Omega$. Therefore, $Aut(A(n,k,k))|_\Omega$ induces a permutation group acting on $\Sigma$, denoted by $(Aut(A(n,k,k))|_\Omega)|_\Sigma$. It is obvious that $(Aut(A(n,k,k))|_\Omega)|_\Sigma\cong S_n$. Consider the natural surjective homomorphism $$\mu:Aut(A(n,k,k))|_\Omega\rightarrow(Aut(A(n,k,k))|_\Omega)|_\Sigma.$$ It follows from Proposition\ \ref{pan2-4} that $ker\mu\cong S_k$. So far, we have seen that $$Aut(A(n,k,k))|_\Omega/S_k\cong S_n.$$ Moreover, $P(S_n)\times Q(S_k)\leq Aut(A(n,k,k))$ and Proposition\ \ref{pan2-2} imply that $$Aut(A(n,k,k))\cong P(S_n)\times Q(S_k).$$ This completes the proof of this lemma.   \qed

Next we consider $A(n,n,n)$ and $A(n,n,2)$. Recall some notions and notations about Cayley graph. Given a group $G$ and a subset $S\subseteq G$ such that $1\not\in S$ and $S^{-1}=S$. The Cayley graph of $G$ with respect to the $S$, denoted by $Cay(G,S)$, is defined to the graph with vertex set $G$ and edge set $\{(g,sg)|g\in G,s\in S\}$. Consider the symmetric group $S_n$. If $T$ is the set of all transpositions of $S_n$, then $Cay(S_n,T)$ is called the complete transposition graph; if $D$ is the set of all derangements of $S_n$, then $Cay(S_n,D)$ is called the complete derangement graph. Ashwin Ganesan \cite{GA}, Yun-Ping Deng and Xiao-Dong Zhang \cite{DZ} characterized the full automorphism groups of $Cay(S_n,T)$ and $Cay(S_n,D)$ respectively. Next, we shall prove that $Cay(S_n,T)\cong A(n,n,2)$ and $Cay(S_n,D)\cong A(n,n,n)$.

\begin{Pro}\label{pan2-6}\normalfont
$Cay(S_n,T)\cong A(n,n,2)$ and $Cay(S_n,D)\cong A(n,n,n)$ for all $n>2$.
\end{Pro}
\demo Defining the mapping $\psi: [n]^{(n)}\rightarrow S_n$ under the rule $\psi([s_1,s_2,...,s_n])=s_1s_2\cdot\cdot\cdot s_n$ for every $n$-tuple $[s_1,s_2,...,s_n]\in[n]^{(n)}$, where $s_1s_2\cdot\cdot\cdot s_n$ is the permutation that maps $i$ to $s_i$ for $i=1,2,...,n$. Clearly, $\psi$ is a bijection between $[n]^{(n)}$ and $S_n$. Next we prove that $\psi$ is not only an isomorphism between $Cay(S_n,T)$ and $A(n,n,2)$, but also an isomorphism between $Cay(S_n,D)$ and $A(n,n,n)$.

Consider $Cay(S_n,T)$ and $A(n,n,2)$. Let $[s_1,s_2,...,s_n]$ and $[t_1,t_2,...,t_n]$ be two vertices of $A(n,n,2)$. We see that $\psi([s_1,s_2,...,s_n])=s_1s_2\cdot\cdot\cdot s_n$ and $\psi([t_1,t_2,...,t_n])=t_1t_2\cdot\cdot\cdot t_n$. By definition of Cayley graph, we know that $\psi([s_1,s_2,...,s_n])$ and $\psi([t_1,t_2,...,t_n])$ are adjacent if and only if $\psi([s_1,s_2,...,s_n])\psi([t_1,t_2,...,t_n])^{-1}\in T$, where $\psi([t_1,t_2,...,t_n])^{-1}$ stands for the inverse of $\psi([t_1,t_2,...,t_n])$. One easily checks that $\psi([s_1,s_2,...,s_n])\psi([t_1,t_2,...,t_n])^{-1}\in T$ if and only if $|\{i|s_i\neq t_i,1\leq i\leq n\}|=2$. Thereby, $Cay(S_n,T)\cong A(n,n,2)$ for all $n>2$. By the same token, we derive $Cay(S_n,D)\cong A(n,n,n)$ for all $n>2$. The proof of this proposition is complete. \qed

\begin{lem}\label{pan2-7}\normalfont
$Aut(A(n,n,2))\cong [P(S_n)\times Q(S_n)]\rtimes Z_2$ and $Aut(A(n,n,n))\cong [P(S_n)\times Q(S_n)]\rtimes Z_2$ for all $n>2$. Here, $Z_2=\langle h\rangle$ and $h$ is the map that $[s_1,s_2,...,s_n]^h=[x_1,x_2,...,x_n]$ where $x_j=i$ if $s_i=j$ for $j=1,2,...,n$.
\end{lem}
\demo Apply Proposition\ \ref{pan2-6} and \cite[Theorem~1.1]{GA} and \cite[Theorem~1.1]{DZ}, we deduce that both the orders of $Aut(A(n,n,2))$ and $Aut(A(n,n,n))$ are $2n!n!$. Moreover, $P(S_n)\times Q(S_n)\leq Aut(A(n,n,2))$ and $P(S_n)\times Q(S_n)\leq Aut(A(n,n,n))$. So it suffices to prove that $h$ is an automorphism of order $2$ that is not in $P(S_n)\times Q(S_n)$. Clearly, $\psi([s_1,s_2,...,s_n]^h)$ is the inverse of $\psi([s_1,s_2,...,s_n])$. By Proposition\ \ref{pan2-6} and \cite[Theorem~1.1]{GA} and \cite[Theorem~1.1]{DZ}, $h$ is an automorphism of $Aut(A(n,n,2))$ and $Aut(A(n,n,n))$ with order $2$. On the other hand, it is easy to see that $\Sigma$ and $\Sigma'$ are two systems of blocks of $P(S_n)\times Q(S_n)|_\Omega$ acting on $\Omega$. However, it is easy to check that $\Sigma$ and $\Sigma'$ are not systems of blocks of $\langle h\rangle|_\Omega$ acting on $\Omega$.
This completes the proof of this lemma. \qed

So far, we have arrived at Theorem\ \ref{pan1-2} by Lemma\ \ref{pan2-5} and Lemma\ \ref{pan2-7}.

\section {Concluding Remarks}

Note that the transposition of $S_n$ has $(n-2)$ fixed points, and derangement of $S_n$ has no fixed point. Thereby, we see that the transposition and derangement have the numbers of fixed points in exactly two extreme cases. However, $Aut(Cay(S_n,T))\cong Aut(Cay(S_n,D))$. So we propose a bold conjecture that is a generalization of the complete transposition graph and derangement graph.
\begin{Con}\label{pan3-1}\normalfont
Let $F_k$ denotes the set of permutations in $S_n$ which have $k$ fixed points with $0\leq k\leq n-2$. Then $Aut(Cay(S_n,F_k))\cong [R(S_n)\rtimes Inn(S_n)]\rtimes Z_2$ for all $n>2$. Here, $R(S_n)$ and $Inn(S_n)$ are the right regular representation and the inner automorphism group of $S_n$ respectively, and $Z_2=\langle h\rangle$ with the mapping $h:g^h=h^{-1},\forall g\in S_n$.
\end{Con}

By the proof of Proposition\ \ref{pan2-6}, we can infer that $Cay(S_n,F_k)\cong A(n,n,n-k)$. In addition, the maximum independent set of $A(n,k,r)$ becomes complicated when $r< k$, however, the Question\ \ref{pan1-1} may be resolved from the perspective of maximum clique. So it seems that continuing to investigate Question\ \ref{pan1-1} is more suitable for confirming Conjecture\ \ref{pan3-1}.

\end{document}